\documentclass{amsart}[12pt]
\usepackage{a4}
\usepackage[T1]{fontenc}
\usepackage[latin1]{inputenc}
\usepackage{hyperref}
\usepackage[dvips]{graphics}
\usepackage[english]{babel}
\usepackage{amsmath}
\usepackage{amssymb}
\usepackage{amsopn}
\usepackage{amsbsy}
\usepackage{amstext}
\usepackage{latexsym}
\usepackage{amsfonts}
\usepackage{array}
\usepackage{epsfig}
\usepackage{mathrsfs}
\title[Any additive basis has only finitely many essential subsets]{A simple proof that any
additive basis has only finitely many essential subsets}
\author[B. Farhi]{Bakir FARHI}
\date{\bf July 22\textsuperscript{nd}, 2008}

\newtheorem{thm}{Theorem}

\newtheorem{lemma}[thm]{Lemma}
\newtheorem{rmk}[thm]{Remark}
\newtheorem{coll}[thm]{Corollary}

\newtheorem{defi}[thm]{Definition}

\def\gcd{{\rm gcd}}
\def\d{{\rm d}}
\def\card{{\rm card}~\!}
\def\mod{{\rm mod}}
\let\epsilon=\varepsilon

\def\EMdash{\leavevmode\hbox to 7.5mm{\vrule height .63ex depth -.59ex
    width 5.4mm\hfill}}

\begin{document}
\maketitle
\begin{center}
Département de Mathématiques, Université du Maine, \\
Avenue Olivier Messiaen, 72085 Le Mans Cedex 9, France. \\
bakir.farhi@gmail.com
\end{center}
\begin{abstract}
Let $A$ be an additive basis. We call ``essential subset'' of $A$
any finite subset $P$ of $A$ such that $A \setminus P$ is not an
additive basis and that $P$ is minimal (for the inclusion order)
to have this property. A recent theorem due to B. Deschamps and
the author states that any additive basis has only finitely many
essential subsets (see ``Essentialit\'e dans les bases additives,
{\it J. Number Theory}, {\bf{123}} (2007), p. 170-192''). The aim
of this note is to give a simple proof of this theorem.
\end{abstract}
\vspace{1mm}

\section{Introduction}

An additive basis (or simply a basis) is a subset $A$ of
$\mathbb{Z}$, having a finite intersection with $\mathbb{Z}^-$ and
for which there exists a natural number $h$ such that any
sufficiently large positive integer can be written as a sum of $h$
elements of $A$. The smaller number $h$ satisfying this property
is called ``the order'' of the basis $A$. Given a basis $A$, an
element $x$ of $A$ is said to be ``essential'' if the set $A
\setminus \{x\}$ is not a basis.

Erd\"os and Graham\cite{eg} proved that an element $x$ of a basis
$A$ is not essential if and only if $\gcd\{a - b ~|~ a , b \in A
\setminus \{x\}\} = 1$. Actually, these two authors proved this
result in the particular case $x = 0 = \min A$, but, as remarked
by Grekos\cite{g}, it suffices to translate $A$ by $(- x)$ to
obtain the generalization. In \cite{nn}, Nash and Nathanson
obtained the following more general result: Let $A$ be an additive
basis and $F$ be a finite subset of $A$. Then the set $A \setminus
F$ is a basis if and only if $\gcd\{a - b ~|~ a , b \in A
\setminus F\} = 1$.

Using the Erd\"os-Graham's characterization, Grekos\cite{g} showed
that the set of the essential elements of a basis $A$ is always
finite and its cardinal can be bounded above in function of the
order of $A$. Recently, Deschamps and the author\cite{df} have
extended the concept of essential element to those of
``essentiality'' and ``essential subset'' which they have defined
as follows:~\vspace{2mm}
\begin{defi}[\cite{df}]\label{d1} Let $A$ be an additive basis. We call ``essentiality'' of $A$ any subset $P$
of $A$ such that $A \setminus P$ is not a basis and that $P$ is
minimal, for the inclusion order, to have this property (so if $Q \subsetneqq P$ then $A \setminus Q$ is a basis). \\
A finite essentiality of $A$ is called an ``essential subset'' of
$A$.
\end{defi}
\noindent{\bf Examples:}
\begin{enumerate}
\item[1)] The set $A = \{6 k ~|~ k \in \mathbb{N}\} \cup \{1 ,
5\}$ is easily seen to be a basis of order $4$. The finite subset
$X = \{1 , 5\}$ of $A$ is an essential subset of $A$, because $A
\setminus X = \{6 k ~|~ k \in \mathbb{N}\}$ is not a basis while
each of the two sets $A \setminus \{1\} = \{6 k ~|~ k \in
\mathbb{N}\} \cup \{5\}$ and $A \setminus \{5\} = \{6 k ~|~ k \in
\mathbb{N}\} \cup \{1\}$ constitutes a basis.
 \item[2)] In the basis $\mathbb{N}$, each of the two complementary subsets $\{2 k ~|~ k \in \mathbb{N}\}$ and
 $\{2 k + 1 ~|~ k \in \mathbb{N}\}$ constitutes an infinite
 essentiality. Indeed, none of those sets is a basis but it
 suffices to add to one of them an element of its complementary to
 obtain a basis.
\end{enumerate}

The number of all essentialities of a basis may be infinite. For
example, we easily verify that for all prime number $p$, the set
$\mathbb{N} \setminus \{p k ~|~ k \in \mathbb{N}\}$ constitutes an
essentiality of the basis $\mathbb{N}$. So, since the set of prime
numbers is infinite then the basis $\mathbb{N}$ contains an
infinitely many essentialities. However, the set of all essential
subsets of a basis is always finite as recently shown by Deschamps
and the author \cite{df} in the following:

\begin{thm}[\cite{df}, Theorem 10]\label{c1}
Any additive basis has only finitely many essential
subsets.~\vspace{2mm}
\end{thm}
In addition, it has been shown in \cite{df} that (contrary to the
set of essential elements) the cardinal of the set of the
essential subsets of an additive basis cannot be bounded above by
a function of the order of the basis alone, but it can be bounded
above in function of another parameter related to the basis.
Below, we give an alternative proof of Theorem \ref{c1}. However,
although our proof is more simple than that of \cite{df}, it does
not permit to bound from above the finite cardinal in question.

\section{A simple proof of Theorem \ref{c1}}

For the following, if $P$ is an essential subset of an additive
basis $A$, we write
$$\d(P) := \gcd\{x - y ~|~ x , y \in A \setminus P\} .$$
Further, if $n$ is a positive integer, we note $\omega(n)$ the
number of its distinct prime factors.

We begin by recalling Lemma 11 of \cite{df}, which constitutes the
main tool of this paper.
\begin{lemma}[Lemma 11 of \cite{df}]\label{l1}
Let $A$ be an additive basis and $P_1$ and $P_2$ be two distinct
essentialities of $A$ such that $P_1 \cup P_2 \neq A$. Then we
have $\d(P_i) \geq 2$ for $i = 1 , 2$ and $\gcd(\d(P_1) , \d(P_2))
= 1$.
\end{lemma}
\noindent{\bf Proof.} Fix $i \in \{1 , 2\}$ and let $x \in P_i$.
Then, because $P_i$ is an essentiality of $A$, the set $(A
\setminus P_i) \cup \{x\}$ is a basis while the set $A \setminus
P_i$ is not a basis. Hence $x$ is an essential element of $(A
\setminus P_i) \cup \{x\}$. This implies (according to the result
of Erd\"os-Graham\cite{eg} and Grekos\cite{g}, cited in §1) that
$\d(P_i) \neq 1$; that is $\d(P_i) \geq 2$ as required.

In order to prove that $\gcd(\d(P_1) , \d(P_2)) = 1$, let us argue
by contradiction. So, assume that there exists $d \geq 2$ such
that $d | \d(P_1)$ and $d | \d(P_2)$. Fix $t \in A \setminus (P_1
\cup P_2)$ and put $B := A \setminus (P_1 \cap P_2)$. For all $x
\in B$, we have $x \not\in P_i$ for some $i \in \{1 , 2\}$, thus
$d | \d(P_i) | (x - t)$, so $x \equiv t ~\mod (d)$. We deduce from
this last fact that $B$ cannot be a basis (because all the
elements of $B$ belong to the same residue class modulo $d \geq
2$). But since $P_1 \cap P_2 \subset P_1$, $P_1 \cap P_2 \subset
P_2$ and $P_1$ and $P_2$ are essentialities of $A$, it follows
that $P_1 \cap P_2 = P_1 = P_2$, which contradicts our hypothesis
that $P_1 \neq P_2$. Hence $\gcd(\d(P_1) , \d(P_2)) = 1$, as
required. The proof is complete.
\penalty-20\null\hfill$\blacksquare$\par\medbreak\noindent
\begin{rmk}\label{r1}
If $P_1$ and $P_2$ are distinct essential subsets of an additive
basis $A$, then the condition $P_1 \cup P_2 \neq A$ of Lemma
\ref{l1} is automatically satisfied (because $A$ is infinite while
$P_1 \cup P_2$ is finite).
\end{rmk}
\begin{coll}\label{c2}
Let $A$ be an additive basis and ${(P_i)}_{i \in I}$ be a nonempty
family of pairwise distinct essential subsets of $A$. Then for all
$(x , y) \in A^2$, with $x \neq y$, the subset of $I$ defined by:
$$J_{x , y} := \{i \in I ~|~ x \not\in P_i ~\text{and}~ y \not\in P_i\}$$
is finite.
\end{coll}
\noindent{\bf Proof.} Let us fix a couple $(x , y)$ of $A^2$ such
that $x \neq y$. From the definition of the set $J_{x , y}$, we
clearly have:
$$\{x , y\} \subset \bigcap_{i \in J_{x , y}}
\left(A \setminus P_i\right) .$$
 This implies that for all $i \in J_{x , y}$, the positive integer $\d(P_i)$ divides the nonzero integer $(x - y)$.
 But since (according to Lemma \ref{l1} and Remark \ref{r1}) the integers $\d(P_i)$ $(i \in J_{x ,
 y})$ are all $\geq 2$ and pairwise coprime, we deduce that their number is at most $\omega(|x - y|)$; so
 $\card J_{x , y} \leq \omega(|x - y|) < + \infty$. The corollary is proved.
 \penalty-20\null\hfill$\blacksquare$\par\medbreak\noindent
{\bf Proof of Theorem \ref{c1}.} Let $A$ be an additive basis and
${(P_i)}_{i \in I}$ be the family of all pairwise distinct
essential subsets of $A$. We have to show that $I$ is finite. If
$\card I \leq 1$ then we are done. Assume for the following that
$\card I \geq 2$ and let us fix $\alpha \in I$. Set for all $x \in
A$:
 $$J_x := \{i \in I ~|~ x \not\in P_i\}$$
and for all $(x , y) \in A^2$:
$$J_{x , y} := \{i \in I ~|~ x \not\in P_i ~\text{and}~ y \not\in P_i\} .$$
Also set $\Lambda$ the finite subset of $A$ defined by:
$$\Lambda := \{x \in P_{\alpha} ~|~ J_x \neq \emptyset\} .$$
This set $\Lambda$ is nonempty (since otherwise we would have
$P_{\alpha} \subset P_i$ $(\forall i \in I)$, which implies
$P_{\alpha} = P_i$ $(\forall i \in I)$, which leads to a
contradiction for any $i \in I$, $i \neq \alpha$, since the
$P_i$'s are pairwise distinct). Now, by the axiom of choice, let
us associate to each element $x \in \Lambda$ (so $J_x \neq
\emptyset$) an element $i(x)$ of $J_x$. We remark that for all $(x
, y) \in A^2$ such that $x \in \Lambda$ and $y \in P_{i(x)}$, we
have $x \neq y$ (because $x \in \Lambda$ implies $i(x) \in J_x$,
that is $x \not\in P_{i(x)}$). It follows from Corollary \ref{c2}
that if a couple $(x , y) \in A^2$ satisfies $x \in \Lambda$ and
$y \in P_{i(x)}$ then the subset $J_{x , y}$ of $I$ is finite.
Consequently, the subset of $I$ defined by:
$$\widetilde{I} := \{\alpha\} \cup \{i(x) ~|~ x \in \Lambda\} \cup \left(\bigcup_{x \in \Lambda , y \in P_{i(x)}}
J_{x , y}\right)$$ is also finite (as a finite union of finite sets).\\
We complete our proof by showing that in fact $I = \widetilde{I}$.
The inclusion $\widetilde{I} \subset I$ is obvious. To show the
second inclusion $I \subset \widetilde{I}$, let us argue by
contradiction; so assume that there exists $i \in I$ such that $i
\not\in \widetilde{I}$. The fact $i \not\in \widetilde{I}$ implies
$i \neq \alpha$ which implies that the two essentialities
$P_{\alpha}$ and $P_i$ are distinct, so $P_{\alpha} \not\subset
P_i$. Thus there exists $x \in P_{\alpha}$ such that $x \not\in
P_i$. Now $x \not\in P_i$ implies $i \in J_x$ which implies $J_x
\neq \emptyset$. Next $x \in P_{\alpha}$ and $J_x \neq \emptyset$
mean that $x \in \Lambda$, hence $i(x) \in \widetilde{I}$. But
since $i \not\in \widetilde{I}$, we certainly have $i \neq i(x)$.
This last fact implies that the two essentialities $P_i$ and
$P_{i(x)}$ are distinct, so $P_{i(x)} \not\subset P_i$. Thus there
exists $y \in P_{i(x)}$ such that $y \not\in P_i$. Finally, the
facts $x \not\in P_i$ and $y \not\in P_i$ imply $i \in J_{x , y}$
which implies (since $x \in \Lambda$ and $y \in P_{i(x)}$) that $i
\in \widetilde{I}$. Contradiction. The proof is complete.
\penalty-20\null\hfill$\blacksquare$\par\medbreak

\end{document}